# Galois Action on Class Groups

Franz Lemmermeyer

May 6, 2002

**Abstract**

It is well known that the Galois group of an extension $L/F$ puts constraints on the structure of the relative ideal class group $\mathrm{Cl}(L/F)$. Explicit results, however, hardly ever go beyond the semisimple abelian case, where $L/F$ is abelian (in general cyclic) and where $(L:F)$ and $\#\mathrm{Cl}(L/F)$ are coprime. Using only basic parts of the theory of group representations, we give a unified approach to these as well as more general results.

It was noticed early on that the action of the Galois group puts constraints on the structure of the ideal class groups of normal extensions; most authors exploited only the action of cyclic subgroups of the Galois groups or restricted their attention to abelian extensions. See e.g. Inaba [10], Yokoyama [24], Iwasawa [12], Smith [19], Cornell & Rosen [5], and, more recently, Komatsu & Nakano [14]. Our aim is to describe a simple and general method that is applicable to arbitrary finite Galois groups.

## 1 Background from Algebraic Number Theory

Let $L$ be a number field, $\mathcal{O}_L$ its ring of integers, $I_L$ its group of fractional ideals $\neq (0)$, and $\mathrm{Cl}(L)$ the ideal class group of $L$ in the usual (wide) sense. In this section we will collect some well known results and techniques that will be generalized subsequently.

For relative extensions $L/F$ of number fields, the relative norm $N_{L/F}: I_L \longrightarrow I_F$ induces a homomorphism $\mathrm{Cl}(L) \longrightarrow \mathrm{Cl}(F)$, which we will also denote by $N_{L/F}$. The kernel $\mathrm{Cl}(L/F)$ of this map is called the *relative class group*. The $p$-Sylow subgroup $\mathrm{Cl}_p(L/F)$ of $\mathrm{Cl}(L/F)$ is the kernel of the restriction of $N_{L/F}$ to the $p$-Sylow subgroup $\mathrm{Cl}_p(L)$ of $\mathrm{Cl}(L)$.

$$\begin{array}{ccc} L & & \mathrm{Cl}(L) \\ | & N_{L/F} \downarrow \uparrow j_{F \to L} & \\ F & & \mathrm{Cl}(F) \end{array}$$

Let $j_{F \to L} : \mathrm{Cl}(F) \longrightarrow \mathrm{Cl}(L)$ denote the transfer of ideal classes induced by mapping ideals $\mathfrak{a} \mathcal{O}_F$ to $\mathfrak{a} \mathcal{O}_L$. Then $N_{L/F} \circ j_{F \to L}(\mathfrak{a}) = \mathfrak{a}^{(L:F)}$. Therefore, the map $\mathrm{Cl}_p(F) \longrightarrow \mathrm{Cl}_p(F) : c \longmapsto c^{(L:F)}$ is an isomorphism for primes $p \nmid (L:F)$, hence $j_{F \to L} : \mathrm{Cl}_p(F) \longrightarrow \mathrm{Cl}_p(L)$ is injective and $N_{L/F} : \mathrm{Cl}_p(L) \longrightarrow \mathrm{Cl}_p(F)$ is surjective; in other words, $j_{F \to L}$ is a section of the exact sequence

$$1 \longrightarrow \mathrm{Cl}_p(L/F) \longrightarrow \mathrm{Cl}_p(L) \xrightarrow{N_{L/F}} \mathrm{Cl}_p(F) \longrightarrow 1.$$

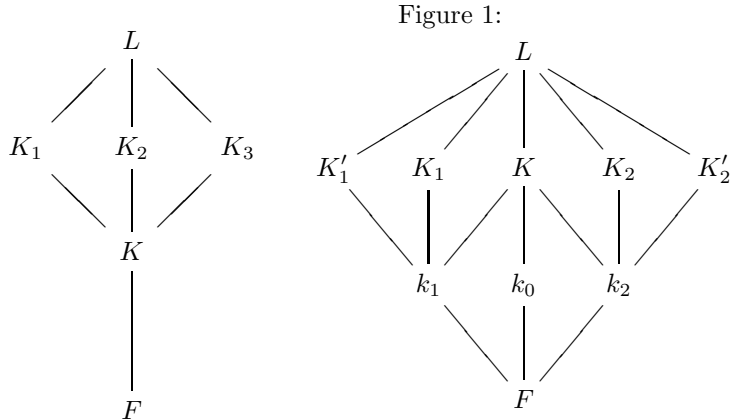

Figure 1:

Thus this sequence splits, and we conclude $\mathrm{Cl}_p(L) \simeq \mathrm{Cl}_p(L/F) \times \mathrm{Cl}_p(F)$.

**Proposition 1.** *Let $L/F$ be an extension of number fields and $p$ a prime not dividing $(L:F)$. Then the transfer of ideal classes $j_{F \to L} : \mathrm{Cl}_p(F) \longrightarrow \mathrm{Cl}_p(L)$ is injective, the relative norm $N_{L/F} : \mathrm{Cl}_p(L) \longrightarrow \mathrm{Cl}_p(F)$ is surjective, and $\mathrm{Cl}_p(L) \simeq \mathrm{Cl}_p(L/F) \times \mathrm{Cl}_p(F)$.*

This result can be used to derive constraints on the class groups of normal extensions; for this type of constraints one needs the presence of conjugate subfields.

Let $V_4 = C_2 \times C_2$ denote Klein's four group, and $C_n$ a cyclic group of order $n$. Moreover, let $A_4$ denote the alternating group of order 12. Hasse diagrams of subfields of extensions with Galois groups $A_4$ and $D_4$ (the dihedral group of order 8) are given in Fig. 1.

**Corollary 2.** *Let $K/F$ be a $V_4$-extension of number fields with quadratic subextensions $k_i$, $i = 0, 1, 2$. Then, for any odd prime $p$,*

$$\mathrm{Cl}_p(K/F) \simeq \mathrm{Cl}_p(k_0/F) \times \mathrm{Cl}_p(k_1/F) \times \mathrm{Cl}_p(k_2/F) \tag{1}$$

*and*

$$\mathrm{Cl}_p(K/k_0) \simeq \mathrm{Cl}_p(k_1/F) \times \mathrm{Cl}_p(k_2/F). \tag{2}$$

*Proof.* First we observe that squaring is an automorphism on every finite abelian group of odd order. Next, let $\sigma_i$ denote the nontrivial automorphism of $k_i/F$. The identity $2 + (1 + \sigma_0 + \sigma_1 + \sigma_2) = (1 + \sigma_0) + (1 + \sigma_1) + (1 + \sigma_2)$ in the group ring $\mathbb{Z}[G]$, $G = V_4$, gives rise to a homomorphism $\mathrm{Cl}_p(L/F) \longrightarrow \mathrm{Cl}_p(k_1/F) \times \mathrm{Cl}_p(k_2/F) \times \mathrm{Cl}_p(k_3/F)$: in fact, write $c \in \mathrm{Cl}_p(L/F)$ in the form $c = d^2$, observe that $d^{(1+\sigma_0+\sigma_1+\sigma_2)} = 1$, and map $c$ to $(d^{(1+\sigma_3)}, d^{(1+\sigma_1)}, d^{(1+\sigma_2)})$. Checking that this is an isomorphism is easy, hence we have (1). Putting together $\mathrm{Cl}_p(K) \simeq \mathrm{Cl}_p(K/k_0) \times \mathrm{Cl}_p(k_0)$ and $\mathrm{Cl}_p(k_0) \simeq \mathrm{Cl}_p(k_0/F) \times \mathrm{Cl}_p(F)$ gives $\mathrm{Cl}_p(K) \simeq \mathrm{Cl}_p(K/k_0) \times \mathrm{Cl}_p(k_0/F) \times \mathrm{Cl}_p(F)$. Comparing this with $\mathrm{Cl}_p(K) \simeq \mathrm{Cl}_p(K/F) \times \mathrm{Cl}_p(F)$ we deduce $\mathrm{Cl}_p(K/F) \simeq \mathrm{Cl}_p(K/k_0) \times \mathrm{Cl}_p(k_0/F)$. Together with (1) this gives (2) as claimed. □

**Corollary 3.** *Let $L/F$ be a dihedral extension of degree $8$, with subextensions denoted as in Fig. 1. Then $\mathrm{Cl}_p(L/K) \simeq \mathrm{Cl}_p(K_1/k_1) \times \mathrm{Cl}_p(K_1/k_1)$ for every odd prime $p$.*

*Proof.* Since $K_1$ and $K_1'$ are conjugate fields (over $F$), they have isomorphic ideal class groups, and by Corollary 2 this implies that $\mathrm{Cl}_p(L/K) \simeq \mathrm{Cl}_p(K_1/k_1) \times \mathrm{Cl}_p(K_1/k_1)$ for every odd prime $p$. This is the algebraic explanation for the fact (following easily from the analytic class number formula) that the order $h_p(L/K)$ of the relative $p$-class group is a square. □

**Corollary 4.** *Let $L/F$ be an $A_4$-extension, with subextensions denoted as in Fig. 1. Then $\mathrm{Cl}_p(L/K) \simeq \mathrm{Cl}_p(K_1/K) \times \mathrm{Cl}_p(K_1/K) \times \mathrm{Cl}_p(K_1/K)$ for odd primes $p$.*

*Proof.* This follows at once from Corollary 2 by observing that the fields $K_i$ are conjugate over $F$. □

## 2 Background from Representation Theory

Here we will review the basic notions of the small part of representation theory that we need. We use Isaacs [11] as our main reference (see also James & Liebeck [13] and Huppert [9]).

Let $F$ be a field and $G$ a finite group; any $F[G]$-module $A$ gives rise to a representation $\Phi : G \longrightarrow \mathrm{Aut}(A)$. Conversely, every homomorphism $\Phi : G \longrightarrow \mathrm{Aut}(A)$ of $G$ into the automorphism group of an $F$-vector space $A$ makes $A$ into an $F[G]$-module, and the submodules of $A$ are exactly the subspaces of $A$ invariant under the action of $G$. Such an $F[G]$-module $A$ is called *irreducible* if $0$ and $A$ are the only submodules of $A$; it is called *completely reducible*, if, for every submodule $B \subseteq A$, there is a submodule $C \subseteq A$ such that $A = B \oplus C$. One of the most basic results of nonmodular representation theory is Maschke's theorem (cf. [11, Thm. 1.9]):

**Theorem 5.** *If $G$ is a finite group, and if $F$ is a field whose characteristic does not divide $\#G$, then all $F[G]$-modules are completely reducible.*

An $F$-representation $\Phi$ of $G$ can be viewed as a homomorphism $G \longrightarrow \mathrm{GL}_n(F)$. In particular, $\Phi(g)$ is a matrix for every $g \in G$, and its trace $\chi(g)$ does not depend on the choice of the basis. We call $\chi$ the *character* associated to $\Phi$. Characters are functions $G \longrightarrow F$ that are constant on conjugacy classes; the number $h$ of conjugacy classes of $G$ is called the *class number* of $G$.

The kernel of a character $\chi$ is defined by $\ker \chi = \{g \in G : \chi(g) = \chi(1)\}$. It is easy to show that $\ker \chi = \ker \Phi$. In particular, $\ker \chi$ is a normal subgroup of $G$. A character or representation with trivial kernel is called *faithful*.

Let $\Phi$ be an $F$-representation of $G$, and suppose that $F \subseteq E$; then we can view $\Phi$ also as an $E$-representation, and we denote this extension by $\Phi^E$. An $F$-representation $\Phi$ of $G$ is called *absolutely irreducible*, if $\Phi^E$ is irreducible for every extension $E$ of $F$. A field $E$ is called a *splitting field* for $G$, if every irreducible $E$-representation of $G$ is absolutely irreducible. Clearly, every algebraically closed field is a splitting field for $G$. For fields with prime characteristic it is not hard to show that finite splitting fields exist for every finite group ([11, Cor. 9.15]):

**Proposition 6.** *Let $G$ be a group with exponent $n$, and let $E$ be an extension of $\mathbb{F}_p$ in which $x^n - 1$ splits into linear factors. Then $E$ is a splitting field for $G$.*

Let $F$ be a field, and suppose that $E \supseteq F$ is a splitting field for the group $G$. Two irreducible characters $\chi$ and $\psi$ of some $E$-representations of $G$ are called *Galois conjugate* over $F$ if $F(\chi) = F(\psi)$, and if there is a $\sigma \in \text{Gal}\,(F(\chi)/F)$ such that $\chi^\sigma = \psi$. This induces an equivalence relation on the set of irreducible $E$-characters. The cardinality of these equivalence classes is given by [11, Lemma 9.17.c)]:

**Proposition 7.** *Let $\Omega$ be the equivalence class of an $E$-character $\chi$ over $F$. Then $\Omega$ has cardinality $(F(\chi) : F)$.*

The notion of faithfulness is compatible with the extension of representations to splitting fields:

**Lemma 8.** *Let $E$ be a splitting field for $G$, $F$ a subfield of $E$, and let $\Phi$ be an irreducible $F$-representation. If $\Phi$ is faithful, then so is each irreducible component of $\Phi^E$.*

*Proof.* Let $\chi$ be the character afforded by $\Phi$. Since $\chi^E$ is the sum of irreducible $E$-characters $\chi_j$ that are Galois conjugate, we have $\chi_j(g) = \chi_j(1)$ for one of these characters if and only if $\chi_j(g) = \chi_j(1)$ for all of them. Thus if one of the components of $\chi^E$ is not faithful, neither is $\chi$. □

Let $E$ be a splitting field for $G$, $F$ a subfield of $E$, and let $\chi$ be an irreducible $E$-character. Then $F(\chi)$ is defined to be the smallest extension of $F$ containing all $\chi(g)$, $g \in G$. Proposition 6 shows that $F(\chi)/F$ is an abelian extension if $F$ is finite. The set of irreducible $E$-characters of $G$ will be denoted by $\text{Irr}_E(G)$. It can be shown that this notion does not depend on the choice of the splitting field $E$ (cf. [11, Lemma 9.13]).

The following theorem describes the behaviour of irreducible representations and their characters in splitting fields (cf. [11, Thm. 9.21]):

**Theorem 9.** *Let $F$ be a field containing $\mathbb{F}_p$, $G$ a finite group such that $p \nmid \#G$, and let $E \supseteq F$ be a splitting field for $G$. If $\Phi$ is an irreducible $F$-representation, then $\Phi^E$ is the direct sum of $E$-irreducible representations, all occurring with multiplicity $1$. Moreover, the characters afforded by the $E$-irreducible components of $\Phi^E$ constitute a Galois conjugacy class over $F$.*

## 3 The Main Theorem

Before we give the main result of this paper, we take the opportunity to recall an almost forgotten result of Grün [7] on the structure of class groups of non-abelian normal extensions:

**Proposition 10.** *Let $L/F$ be a normal extension of number fields, and let $K$ denote the maximal subextension of $L/F$ that is abelian over $F$.*

  i) *If $\mathrm{Cl}(L/K)$ is cyclic, then $h(L/K) \mid (L:K)$;*

  ii) *If $\mathrm{Cl}(L)$ is cyclic, then $h(L) \mid (L:K)e$, where $e$ denotes the exponent of $N_{L/K}\mathrm{Cl}(L)$. Observe that $e \mid h(K)$, and that $e = 1$ if $L$ contains the Hilbert class field $K^1$ of $K$.*

*Similar results hold for the p-Sylow subgroups.*

*Proof.* Let $C$ be a cyclic group of order $h$ on which $\Gamma = \mathrm{Gal}\,(L/F)$ acts. This is by definition equivalent to the existence of a homomorphism $\Phi : \Gamma \longrightarrow \mathrm{Aut}(C) \simeq \mathbb{Z}/(h-1)\mathbb{Z}$. Since $\mathrm{im}\,\Phi$ is abelian, $\Gamma' \subseteq \ker \Phi$, hence $\Gamma'$ acts trivially on $C$. Now $\Gamma'$ corresponds to the field $K$ via Galois theory, and we find $N_{L/K}c = c^{(L:K)}$ for all $c \in C$. Putting $C = \mathrm{Cl}(L/K)$, we see at once that $(L:K)$ annihilates $\mathrm{Cl}(L/K)$. If we denote the exponent of $N_{L/K}\mathrm{Cl}(L)$ by $e$ and put $C = \mathrm{Cl}(L)$, then we find in a similar way that $(L:K)e$ annihilates $\mathrm{Cl}(L)$. □

We now come to the main theorem:

**Theorem 11.** *Let $L/F$ be a normal extension of number fields with Galois group $\Gamma = \mathrm{Gal}\,(L/F)$. Let $p$ be a prime not dividing $\#\Gamma$, and assume that $\mathrm{Cl}_p(K/F) = 1$ for all normal subextensions $K/F$ with $F \subseteq K \subsetneq L$. Let $\Phi$ denote the representation $\Phi : \Gamma \longrightarrow \mathrm{Aut}(C)$ induced by the action of $\Gamma$ on $C = \mathrm{Cl}(L/F)/\mathrm{Cl}(L/F)^p$; then $\Phi$ is faithful.*

*If the degrees of all irreducible faithful $\mathbb{F}_p$-characters $\chi$ of $\Gamma$ are divisible by $f$, then $\mathrm{rank}\,\mathrm{Cl}_p(L) \equiv 0 \bmod f$. If, in addition, $(\mathbb{F}_p(\chi) : \mathbb{F}_p) = r$ for all these $\chi$, then $\mathrm{rank}\,\mathrm{Cl}_p(L) \equiv 0 \bmod rf$.*

*Proof.* Clearly $G = \ker \Phi$ is a normal subgroup of $\Gamma$ which acts trivially on $C$. Let $K$ be the subextension of $L/F$ corresponding to $G$ by Galois theory. Then, for every $c \in C$ we have $N_{L/K}c = c^{(L:K)}$, because the automorphisms of $L/K$ leave $c$ invariant. On the other hand $N_{L/K}c$ is an element of $\mathrm{Cl}(K)/\mathrm{Cl}(K)^p$, where $K/F$ is a normal subextension of $L$. Since $\mathrm{Cl}_p(K/F) = 1$ by assumption, we conclude that $c = 1$. Thus $c$ is killed by $(L:K)$ and by $p$; since $p \nmid (L:F)$, we find that $(L:K) = 1$, i.e. $G = 1$, and $\Phi$ is faithful.

If all faithful $\mathbb{F}_p$-characters $\chi$ of $\Gamma$ have degrees divisible by $f$, then clearly $\operatorname{rank} \operatorname{Cl}_p(L) \equiv 0 \bmod f$ since every character is the sum of irreducible characters. If, in addition, $(\mathbb{F}_p(\chi) : \mathbb{F}_p) = r$ for all these $\chi$, and if $\Phi$ denotes an irreducible faithful representation, then Theorem 9 shows that $\Phi^E$ is the direct sum of $r$ representaions of degree $f$, which in turn implies that $\Phi$ has degree $rf$. □

## 4 Applications

We first derive a couple of corollaries from Theorem 11:

**Corollary 12.** *Let $L/F$ be a normal $\ell$-extension with $\Gamma = \operatorname{Gal}(L/F)$, and assume that the center $Z(\Gamma)$ is not cyclic. If $p \neq \ell$ is a prime dividing $h(L)$, then $p$ divides the class number of some normal subextension $K/F$ of $L/F$.*

*Proof.* If not, then there exists a faithful irreducible representation of $\Gamma$. According to [9, V, §5, Ex. 15], an $\ell$-group $\Gamma$ possesses faithful irreducible representations if and only if $Z(\Gamma)$ is cyclic. This contradicts our assumption. □

**Corollary 13.** *Let $L/F$ be a normal $\ell$-extension with nonabelian Galois group $\Gamma = \operatorname{Gal}(L/F)$. If the prime $p \neq \ell$ does not divide the class number of any normal subextension $K/F$ of $L/F$, then $\operatorname{rank} \operatorname{Cl}_p(L) \equiv 0 \bmod \ell$.*

*Proof.* It is sufficient to show that every irreducible faithful character of $\Gamma$ has degree divisible by $\ell$. To this end it is sufficient to show that linear characters cannot be faithful (recall that the degree of a character over a splitting field divides the group order). So assume that $\chi$ is linear and faithful; linear characters are lifts from characters of $\Gamma/\Gamma'$, and this group contains a subgroup of type $(\ell, \ell)$ since $\Gamma$ is nonabelian. In particular, $\Gamma/\Gamma'$ has noncyclic center and, a fortiori, no faithful characters. □

The theorem below is well known in the case where $\Gamma$ is a cyclic group; the other examples serve only to illustrate the method.

**Theorem 14.** *Let $L/F$ be a normal extension with Galois group $\Gamma$, and let $p \nmid \#\Gamma$ be a prime. Suppose that $\operatorname{Cl}_p(K/F) = 1$ for all normal subextensions $K/F$ of $L/F$ with $K \neq L$. Then $\operatorname{rank} \operatorname{Cl}_p(L/F)$ is divisible by $r$, where $r$ is given in the table below:*

| $\Gamma$ | $\#\Gamma$ | $r$ | $f > 0$ *minimal with* |
|---|---:|---:|---|
| $C_n$ | $n$ | $f$ | $p^f \equiv 1 \bmod n$ |
| $D_n$, $n \geq 3$ | $2n$ | $2f$ | $p^f \equiv \pm 1 \bmod n$ |
| $H_{4n}$, $n \geq 2$ | $2n$ | $2f$ | $p^f \equiv \pm 1 \bmod 2n$ |
| $C_4 \curlyvee H_8$ | $16$ | $2f$ | $p^f \equiv 1 \bmod 4$ |
| $A_4$ | $12$ | $3$ | |

*Proof.* This follows immediately from Theorem 11 and character tables for the corresponding groups.

For $\Gamma = C_n = \langle a \rangle$ (see [13, p. 82]), all irreducible characters $\chi_j$ ($0 \leq j \leq n-1$) are linear and given by $\chi_j(a^k) = \zeta_n^{kj}$; the character $\chi_j$ is faithful if and only if $\gcd(j, n) = 1$, and then $\mathbb{F}_p(\chi_j) = \mathbb{F}_p(\zeta_n)$. By the decomposition law for cyclotomic extensions we know that $(\mathbb{F}_p(\zeta_n) : \mathbb{F}_p)$ is just the order of $p$ mod $n$, and this implies the claim if $\Gamma = C_n$.

For $\Gamma = D_n$ ([13, p. 182, 183]), the faithful irreducible characters $\chi$ have degree 2 and satisfy $\mathbb{F}_p(\chi) = \mathbb{F}_p(\zeta_n + \zeta_n^{-1})$. Thus $(\mathbb{F}_p(\chi) : \mathbb{F}_p) = f$ where $f > 0$ is the minimal integer such that $p^f \equiv \pm 1 \mod n$. Since each irreducible faithful $\Phi$ consists of $f$ irreducible faithful characters $\chi$ with degree 2, the claim follows.

For generalized quaternion groups $\Gamma = H_{4n}$ ([13, p. 385]), the faithful irreducible characters $\chi$ have degree 2 and satisfy $\mathbb{F}_p(\chi) = \mathbb{F}_p(\zeta_{2n} + \zeta_{2n}^{-1})$. Thus $(\mathbb{F}_p(\chi) : \mathbb{F}_p) = f$ where $f > 0$ is the minimal integer such that $p^f \equiv \pm 1 \mod 2n$.

Finally, consider the group $\Gamma = A_4$ with four conjugacy classes $C_1$, $C_2$, $C_3$, $C_4$ and the character table (see [13, p. 180, 181]; $\rho$ denotes a primitive cube root of unity)

|          | $C_1$ | $C_2$ | $C_3$   | $C_4$   |
|----------|-------|-------|---------|---------|
| $\chi_1$ | 1     | 1     | 1       | 1       |
| $\chi_2$ | 1     | 1     | $\rho$  | $\rho^2$|
| $\chi_3$ | 1     | 1     | $\rho^2$| $\rho$  |
| $\chi_4$ | 3     | $-1$  | 0       | 0       |

Clearly $\chi_4$ is the only irreducible faithful character, which implies the claims. Note that, in this case, the presence of conjugate subfields is responsible for the divisibility of the $p$-rank of $\mathrm{Cl}(L/F)$ by 3 (compare Corollary 4). □

**Remark.** By replacing $A$ with $A = C^{p^m}/C^{p^{m+1}}$, we can give similar results on the $p^m$-rank of $\mathrm{Cl}(L/F)$.

We also can obtain results for groups with irreducible faithful characters of different degrees by adding assumptions about the triviality of $\mathrm{Cl}_p(K/F)$ for certain nonnormal subextensions.

## 4.1 Cyclic Extensions

The following proposition is well known (cf. [10], [6]):

**Proposition 15.** *Let $p$ be an odd prime, and suppose that $L/K$ is a cyclic extension of degree $p$. If $\mathrm{Cl}_p(L/K)$ is cyclic, then $\mathrm{Cl}_p(L/K)$ is trivial or $\simeq \mathbb{Z}/p\mathbb{Z}$.*

In fact there are even stronger results describing $\mathrm{Cl}_p(L/K)$ as a $\mathrm{Gal}\,(L/K)$-module. The following is a slight generalization for $p = 2$:

**Proposition 16.** *Let $L/F$ be a cyclic quartic extension with Galois group $\Gamma = \mathrm{Gal}\,(L/F) = \langle \sigma \rangle$, and let $K$ be its quadratic subextension, i.e. the fixed field of $\langle \sigma^2 \rangle$. If $C$ is a cyclic subgroup of $\mathrm{Cl}_2(L/K)$ and a $\Gamma$-module, then $\#C \mid 2$.*

*Proof.* Assume that $C$ has an ideal class $c$ of order 4. Then the order of $c^\sigma$ is also 4, hence we must have $c^\sigma = c$ or $c^\sigma = c^3$. In both cases we get $c = c^{\sigma^2}$, hence $1 = N_{L/K} c = c^{1+\sigma^2} = c^2$, contradicting our assumption. Therefore, $C$ is elementary abelian, and since it is cyclic by assumption, our claim follows. □

As an application, we note

**Corollary 17.** *Let $L/\mathbb{Q}$ be a cyclic quartic extension with quadratic subfield $K$. If $\mathrm{Cl}_2(L/K) \simeq (2^\alpha, 2^\beta, \ldots)$, where $\alpha \geq \beta \geq \ldots$, then $\alpha - \beta \in \{0, 1\}$.*

*Proof.* Let $c \in \mathrm{Cl}_2(L/K)$ be an ideal class of order $2^\alpha$. Then $c^{2^\beta}$ generates a subgroup $C$ of $\mathrm{Cl}_2(L/K)$. Since $c^\sigma = c^j \cdot d$ for some odd integer $j$ and an ideal class $d$ of order dividing $2^\beta$, we conclude that $C$ is a $\mathrm{Gal}\,(L/\mathbb{Q})$-module. Since $C$ is cyclic, Proposition 16 implies that $\#C \leq 2$, and this proves our claim. □

This strengthens results of [2, 3] on the structure of 2-class groups of certain quartic cyclic number fields.

## 4.2 Quaternion Extensions

Next assume that $L/F$ is a normal extension with Galois group $\mathrm{Gal}\,(L/F) \simeq H_8$, the quaternion group of order 8. Let $K$ be the unique quartic subextension of $L/F$. Louboutin [18] has computed some relative class numbers $h(L/K)$ in the special case where $F = \mathbb{Q}$ and $L$ is a CM-field: it is known that, for primes $q \equiv 3 \bmod 8$, the real bicyclic biquadratic number field $K = \mathbb{Q}(\sqrt{2}, \sqrt{q}\,)$ admits a unique totally complex quadratic extension $L_q/K$ with the following properties:

a) $L_q/\mathbb{Q}$ is normal, and $\mathrm{Gal}\,(L_q/\mathbb{Q}) \simeq H_8$;

b) $L_q/K$ is unramified outside $\{\infty, 2, q\}$.

Here are his results:

| $q$ | 3 | 11 | 19 | 43 | 59 | 67 | 83 |
|---|---|---|---|---|---|---|---|
| $h(L_q)$ | 2 | $2 \cdot 3^2$ | $2 \cdot 7^2$ | $2 \cdot 3^4$ | $2 \cdot 3^2 \cdot 7^2$ | $2 \cdot 3^6$ | $2 \cdot 5^4$ |

Louboutin observed that the class numbers $h_p(L/K)$, $p$ odd, were all squares. In contrast to the cases $G = A_4$ or $G = D_4$ (cf. Corollaries 3, 4), however, this cannot be explained by the presence of subfields.

For primes $p \equiv 3 \bmod 4$, we see at once that the $p$-part of $\mathrm{Cl}(L/K)$ has even rank, because $L$ is a cyclic quartic extension over each of its three quadratic subextensions, and the action of any of the cyclic groups proves our claim. For primes $p \equiv 1 \bmod 4$, Proposition 10 shows that $\mathrm{Cl}_p(L/K)$ cannot be cyclic, but it does not allow us to deduce that its rank is even: this only follows from the results in the table of Theorem 14. In this special case of quaternion extensions of degree 8 this has already been noticed by Tate (cf. [20, p. 51, Lemme 2.2]).

## 4.3 Dihedral Extensions

We conjecture that the result in the table of Theorem 14 pertaining to $D_n$ also holds for $p = 2$, that is: if the relative class number of the quadratic subextension of $L/F$ is odd, then the 2-rank of $\mathrm{Cl}(L/F)$ is divisible by $2f$, where $f > 0$ is minimal with $2^f \equiv \pm 1 \bmod n$. This is certainly true if $2^f \equiv -1 \bmod n$, because in this case the claim follows from the action of the cyclic subgroup of $D_n$ alone. Ken Yamamura [23] gives ad hoc proofs for some special cases when $2^f \equiv 1 \bmod n$ and $f$ is odd. The general case requires studying modular representations, and I hope to return to this question at another occasion.

Amazingly, our results on dihedral extensions can be used to show that there are constraints on the structure of class groups of non-normal extensions:

**Corollary 18.** *Let $K/F$ be an extension of odd degree $n$ such that its normal closure $L/F$ has Galois group $G = D_n$, and let $r_p$ denote the rank of $\mathrm{Cl}_p(K/F)$. Then $r_p \equiv 0 \bmod f$ for every prime $p \nmid 2n$, where $f$ is the smallest positive integer such that $p^f \equiv \pm 1 \bmod n$. This result is also valid for $p = 2$ if $2^f \equiv -1 \bmod n$.*

*Proof.* Put $A = \mathrm{Cl}_p(K/F)$ and let $k$ be the quadratic subextension of $L/F$; the map $j_{K \to L} : \mathrm{Cl}_p(K/F) \longrightarrow \mathrm{Cl}_p(L/k)$ is known to be injective (for $p \nmid 2n$ this follows from Proposition 1; for $p = 2$ it is due to J.-F. Jaulent, cf. [4, Thm. 7.8]), so we can regard $A$ as a subgroup in $\mathrm{Cl}_p(L/k)$. Let $C$ be the $G$-module generated by $j(A)$. According to [8, Lemma 1], we have $C = AA^\sigma$, where $\sigma$ generates $\mathrm{Gal}\,(L/k)$. Now we claim that $C = A \oplus A^\sigma$ as abelian groups (observe that $A$ and $A^\sigma$ in general are not $G$-modules). To this end we have to show that $A \cap A^\sigma = \{1\}$. Let $a \in A \cap A^\sigma$; from $a = b^\sigma$ for some $b \in A$ we get $a = a^\tau = b^{\sigma\tau} = b^{\tau\sigma^{-1}} = b^{\sigma^{-1}} = (b^\sigma)^{\sigma^{-2}} = a^{\sigma^{-2}}$, hence $a = a^\sigma$. But then $1 = N_{K/F}a = a^n$ contradicting our assumption that $a \in \mathrm{Cl}_p(K/F)$ and $p \nmid n$.

Thus $\mathrm{rank}\,C = \mathrm{rank}\,A + \mathrm{rank}\,A^\sigma = 2 \cdot \mathrm{rank}\,A$. On the other hand, we have $\mathrm{rank}\,C \equiv 0 \bmod 2f$ by the results of the table in Theorem 14. $\square$

**Examples:** We have used PARI to compute class numbers of some dihedral extensions. Kondo [15] describes a family of quintic dihedral extensions whose normal closure is unramified over its quadratic subfield: it is given by $f(x) = x^5 + (a-3)x^4 + (b-a+3)x^3 + (a^2 - a - 1 - 2b)x^2 + bx + a$, where $a$ and $b$ are integers. Putting $a = 1$, we find a subfamily which can be viewed as a series of 'simplest quintic dihedral fields', since these have a parametrized system of independent units.

## 4.4 A Nonabelian Group of Order 16

The group $\Gamma = C_4 \curlyvee H_8 = 16.08$ (see the tables by Thomas & Wood [22]) is a 2-group of order 16. Let $\sigma$ denote the generator of its commutator group $\Gamma'$ of order 2. $\Gamma$ has 8 linear characters, and all of them satisfy $\chi(\sigma) = 1$. $\Gamma$ has two conjugate faithful characters $\psi$ and $\psi'$ of dimension 2, and we have $F(\psi) = F(\sqrt{-1}\,)$. These facts show

Table 1:

| $b$ | disc $K$ | $h$ | $b$ | disc $K$ | $h$ |
|---|---:|---:|---|---:|---:|
| 0 | 2209 | [1] | 10 | 50225569 | [19] |
| 1 | 10609 | [1] | 11 | 81414529 | [4, 4] |
| 2 | 57121 | [1] | 12 | 127215841 | [11] |
| 3 | 229441 | [1] | 13 | 192626641 | [4, 4] |
| 4 | 717409 | [5] | 14 | 283821409 | [19] |
| 5 | 1868689 | [2, 2] | 15 | 408322849 | [10, 2] |
| 6 | 4255969 | [2, 2] | 16 | 575184289 | [61] |
| 7 | 8755681 | [3, 3] | 17 | 795183601 | [4, 4] |
| 8 | 16638241 | [2, 2] | 18 | 450241 | [1] |
| 9 | 29669809 | [11] | 19 | 1447574209 | [20, 4] |

**Proposition 19.** *Let $L/k$ be a normal extension of number fields such that $\Gamma = \mathrm{Gal}\,(L/k) \simeq C_4 \curlyvee H_8$, and let $K$ be the fixed field of $\Gamma' = \langle \sigma \rangle$. Then*

$$\mathrm{rank}\,\mathrm{Cl}_p(L/K) \equiv \begin{cases} 0 \bmod 2, & \text{if } p \equiv 1 \bmod 4, \\ 0 \bmod 4, & \text{if } p \equiv 3 \bmod 4. \end{cases}$$

We can apply this result to find conditions on class groups of extensions of larger degree:

**Corollary 20.** *Let $k$ be a quadratic number field, and let $L/k$ be an unramified extension such that $\Gamma = \mathrm{Gal}\,(L/k) \simeq 32.40$, and assume that the two quaternion extensions of $k$ inside $L$ are normal over $\mathbb{Q}$. Let $K$ be the fixed field of the commutator subgroup $\Gamma' \simeq (2,2)$. Then*

$$\mathrm{rank}\,\mathrm{Cl}_p(L/K) \equiv \begin{cases} 0 \bmod 2, & \text{if } p \equiv 1 \bmod 4, \\ 0 \bmod 4, & \text{if } p \equiv 3 \bmod 4. \end{cases}$$

*Proof.* $\Gamma' \simeq (2,2)$, hence $L/K$ is a $V_4$-extension. Let $L_1, L_2, L_3$ denote the three quadratic subextensions of $L/K$. Then $\mathrm{Gal}\,(L_1/k) \simeq \mathrm{Gal}\,(L_2/k) \simeq D_4 \curlyvee C_4$ and $\mathrm{Gal}\,(L_3/k) \simeq C_2 \times H_8$. Moreover we know from Section 1 that

$$\mathrm{Cl}_p(L/K) \simeq \mathrm{Cl}_p(L_1/K) \times \mathrm{Cl}_p(L_2/K) \times \mathrm{Cl}_p(L_3/K).$$

Applying Proposition 19 to $L_j/k$ for $j = 1, 2$ we see that $\mathrm{rank}\,\mathrm{Cl}_p(L_j/K)$ is even, and divisible by 4 if $p \equiv 3 \bmod 4$. It is therefore sufficient to prove the same for $\mathrm{Cl}_p(L_3/K)$. Now $G = \mathrm{Gal}\,(L_3/K) \simeq C_2 \times H_8$. $G$ has three subgroups $G_i$ of index 2, and these subgroups correspond to three extensions $K$, $K_2$ and $K_3$, which in turn are quadratic extensions of a field $F$ which corresponds to the intersection of the $G_i$. Since $L_3/F$ is a $V_4$-extension, we find

$$\mathrm{Cl}_p(L_3/K) \simeq \mathrm{Cl}_p(K_2/F) \times \mathrm{Cl}_p(K_3/F).$$

It is easily checked that $\mathrm{Gal}\,(K/k) \simeq (2,2,2)$ and $\mathrm{Gal}\,(K_j/k) \simeq H_8$ for $j = 1, 2$. Since the $K_j$ are normal over $\mathbb{Q}$ by assumption, we have $\mathrm{Gal}\,(K_j/\mathbb{Q}) \simeq D_4 \curlyvee C_4$ (cf. [16, 17]). Applying Proposition 19 to $K_j/\mathbb{Q}$ completes our proof. □

As an application we give a proof for the fact that the class field tower of $k = \mathbb{Q}(\sqrt{-105}\,)$ terminates with $k^2$ which differs from the one in [16]. We know that $\mathrm{Cl}(k) \simeq (2,2,2)$, $\mathrm{Cl}(k^1) \simeq (2,2)$, and $\Gamma = \mathrm{Gal}\,(k^2/k) \simeq 32.40$. We also know that $h(k^2)$ is odd from Proposition 21. If $\mathrm{Cl}(k^2)$ were non-trivial, Corollary 20 would show that $h(k^2) \geq 25$. This contradicts the unconditional Odlyzko bounds, which show that $h(k^2) \leq 10$.

The proposition referred to is

**Proposition 21.** *If $k$ is a number field such that $\mathrm{Cl}_2(k_2^1) \simeq (2,2)$, then $k_2^3 = k_2^2$.*

*Proof.* Suppose that $k_2^3 \neq k_2^2$ and let $G = \mathrm{Gal}\,(k_2^3/k)$. Then a result of Taussky [21] shows that $G' \simeq \mathrm{Gal}\,(k_2^3/k_2^1)$ is dihedral, semidihedral or quaternionic, and all these groups have cyclic centers. Burnside [1] has shown that $p$-groups $G$ with cyclic $Z(G')$ have cyclic $G'$. But if $k_2^3/k_2^1$ is cyclic, we must have $k_2^3 = k_2^2$ in contradiction to our assumption. □

### Acknowledgement

I would like to thank A. Brandis and R. Schoof for some helpful discussions, and T. Kondo for sending me [15].